\documentclass{article}

\usepackage{amsmath}
\usepackage{easy}
\usepackage[dvips]{graphics}

\newcommand{\bold}[1]{ \mbox{\boldmath{$#1$}} }

\begin{document}

\title{Observed Range Maximum Likelihood Estimation}
\author{Plamen Markov\\ 
pvmarkov@gmail.com}
\maketitle

\begin{abstract}
The idea of maximizing the likelihood of the observed range for a set of jointly realized counts has been employed in a variety of contexts. The applicability of the MLE introduced in \cite{r1} has been extended to the general case of a multivariate sample containing interval censored outcomes. In addition, a kernel density estimator and a related score function have been proposed leading to the construction of a modified Nadaraya-Watson regression estimator. Finally, the author has treated the problems of  estimating the parameters of a mutinomial distribution and the analysis of contingency tables in the presence of censoring.
\end{abstract}
\section{Summary of previous work}

Let $X_1, X_2, \ldots, X_N$ be i.i.d. real valued random variables with distribution function $F$ and corresponding realized values $x_1, x_2, \ldots, x_N$. In the remainder of the paper we assume that $n \in \{1, 2, \ldots, N\}$. The realized value $x_n$ of the random variable $X_n$ is either an exact observation or censored into an interval $(t_n, t_{2n}]$. We allow for the possibility that $t_{2n}=\infty$ and adopt the convention that $(t_n, t_{2n}]$ is to be interpreted as $(t_n, \infty)$ in that special case.

For a given element $\tau \in \mbox{dom}(F)$ we define $d_\tau$ as the number of sample values observed to be less than or equal to $\tau$ and $a_\tau$ as the number of sample values observed to be greater than $\tau$. The count $u_\tau$ represents the number of censored sample values with censoring intervals that capture $\tau$. For example, a censored value $x_n$ is included in the count $d_\tau$ iff $t_{2n}\leq\tau$ and in the count $a_\tau$ iff $t_n\leq\tau$. From these definitions immediately follows that for any $\tau \in \mbox{dom}(F)$ we have that $d_\tau+a_\tau+u_\tau = N$.

Let $k_\tau$ be the actual number of sample values not exceeding $\tau$. Due to the presence of the censoring mechanism the value of $k_\tau$ is only observed to satisfy $d_\tau \leq k_\tau \leq d_\tau+u_\tau$; we label the latter event as $E$. The likelihood of $E$ is given by
\[L(F(\tau); E) = \sum^{d_\tau + u_\tau}_{k_\tau = d_\tau}  {N \choose k_\tau}\left[F(\tau)\right]^{k_\tau} \left[1-F(\tau)\right]^{N-k_\tau}\]
Let us define the function $\hat{F}: \mbox{dom}(F) \rightarrow [0, 1]$ as the value of $p$ that maximizes
\[L(p; E) = \sum^{d_\tau+u_\tau}_{k=d_\tau}  {N \choose k}(p)^k (1-p)^{N-k}\]
subject to the constraint $0 \leq p \leq 1$. The value of $\hat{F}(\tau)$ has been derived to be
\[ \hat{F}(\tau) = \left\{
	\begin{array}{lll}
		0 & \mbox{if} & d_{\tau}=0 \mbox{ and } a_{\tau}\geq 1 \\ 
		\\
		1 & \mbox{if} & a_{\tau}=0 \mbox{ and } d_{\tau}\geq 1 \\
		\\
		1/2 & \mbox{if} & u_{\tau}=N \\
		\\		\left(1+\sqrt[u_{\tau}+1]{\frac{a_{\tau}(a_{\tau}+1)\ldots(a_{\tau}+u_{\tau})}{d_{\tau}(d_{\tau}+1)\ldots(d_{\tau}+u_{\tau})}}\,\right)^{-1}&\mbox{o.w.}
	\end{array}
\right. \]
Furthermore, the function $\hat{F}$ can be used as an estimator for $F$ since it is a non-decreasing function over $\mbox{dom}(F)$.
\section{Multivariate extension}

In this section the definition of the estimator $\hat{F}$ has been extended to the general case of a sample of $M$-variate observations. Let $\bold{X}_1, \bold{X}_2, \ldots, \bold{X}_N$ be i.i.d. $M$-vectors with distribution function $F$ and the matrix $\bold{D}$ be defined as
\[ \bold{D} = 
\left( \begin{array}{c}
\bold{X}_1\\
\bold{X}_2\\
\vdots\\
\bold{X}_N
\end{array} \right) =
\left( \begin{array}{cccc}
X_{11} & X_{12} & \ldots & X_{1M}\\
X_{21} & X_{22} & \ldots & X_{2M}\\
\vdots & \vdots & & \vdots\\
X_{N1} & X_{N2} & \ldots & X_{NM}
\end{array} \right)
\]
For the rest of the paper we have assumed that all observations $X_{nm}$ are censored into corresponding intervals $(L_{nm}, R_{nm}]$ since the treatment of a dataset $\bold{D}$ containing exact in addition to censored observations does not provide any new mathematical insight.

We also adopt the convention that unless explicitly stated otherwise, an index represented by a small letter ranges between $1$ and the value of the corresponding capital letter inclusive. For example, $m \in \{1, 2, \ldots, M\}$. Furthermore, a random quantity will be always designated by a capital letter and the corresponding small letter will be reserved for its realization. For example, $\bold{x}_n$ is the realization of the random vector $\bold{X}_n$.

Let $\tilde{X}_{i_m}^{(m)}$ be the value of the $i_m$-th biggest element, $i_m \in \{1, 2, \ldots, I_m\}$, of the set
\[ \left\{ L_{1m},\, L_{1m}, \dots, L_{Nm} \right\} \cup \left\{ R_{1m},\, R_{1m}, \dots, R_{Nm} \right\} \]
and the set $G^{(m)}$ be defined as
\[ G^{(m)} = \left\{ \tilde{X}_1^{(m)},\, \tilde{X}_2^{(m)}, \ldots, \tilde{X}_{I_m}^{(m)} \right\} \]
Consequently, the elements of $G^{(m)}$ are all distinct and such that
\[ \tilde{X}_1^{(m)} < \tilde{X}_2^{(m)} < \ldots < \tilde{X}_{I_m}^{(m)} \]
Let us also define the grid $G$ as $G = G^{(1)} \times G^{(2)} \times \ldots \times G^{(M)}$. Our goal will be to estimate $F$ over $G$.

Let $\bold{x} = (x_1, x_2, \ldots, x_M) \in R^M$ and $\bold{x}' = (x'_1, x'_2, \ldots, x'_M) \in R^M$. We will write $\bold{x} < \bold{x}'$ iff $x_m < x'_m$. The expressions $\bold{x} > \bold{x}'$, $\bold{x} \leq \bold{x}'$ and $\bold{x} \geq \bold{x}'$ are defined analogously. Let $\bold{L}_n = (L_{n1}, L_{n2}, \ldots, L_{nM})$ and $\bold{R}_n = (R_{n1}, R_{n2}, \ldots, R_{nM})$. By analogy with the $1$-dimensional case we define $d(\bold{x})$ as the count of observations $\bold{X}_n$ such that $\bold{R}_n \leq \bold{x}$ and $u(\bold{x})$ as the count of observations satisfying $\bold{L}_n < \bold{x} < \bold{R}_n$. It is important to point out that the count $a(\bold{x}) = N - d(\bold{x}) - u(\bold{x})$ is not the number of observations such that $\bold{R}_n < \bold{x}$. Finally, let $k(\bold{x})$ be the realized value of the actual count of observations such that $\bold{R}_n \leq \bold{x}$ and $E$ designate the event $d(\bold{x})\leq k(\bold{x})\leq d(\bold{x})+u(\bold{x})$.

Now we can estimate $F(\bold{x})$ by the value of the variable $p$ that maximizes the function
\[L(p; E) = \sum^{d(\bold{x})+u(\bold{x})}_{k=d(\bold{x})}  {N \choose k}(p)^k (1-p)^{N-k}\]
subject to the constraint $0 \leq p \leq 1$. Consequently the estimator $\hat{F}$ of the unknown distribution function $F$ is given by
\[ \hat{F}(\bold{x}) = \left\{
	\begin{array}{lll}
		0 & \mbox{if} & d(\bold{x})=0 \mbox{ and } a(\bold{x})\geq 1 \\ 
		\\
		1 & \mbox{if} & a(\bold{x})=0 \mbox{ and } d(\bold{x})\geq 1 \\
		\\
		1/2 & \mbox{if} & u(\bold{x})=N \\
		\\		\left(1+\sqrt[u\left(\bold{x}\right)+1]{\frac{a(\bold{x})(a(\bold{x})+1)\ldots(d(\bold{x})+u(\bold{x}))}{d(\bold{x})(d(\bold{x})+1)\ldots(d(\bold{x})+u(\bold{x}))}}\,\,\right)^{-1}&\mbox{o.w.}
	\end{array}
\right. \]

We briefly consider once again a sample of univariate observations $X_1, \ldots, X_N$ with $X_n$ censored into an interval $(L_n, R_n]$ and assume that the random vectors $(X_n,\,L_n,\,R_n)$ are all i.i.d. according to some cdf $F_{XLR}$. The latter function provides a quantitative descsription of the censoring mechanism at play. By setting $\bold{X}_n = (X_n,\,L_n,\,R_n)$ and employing the estimation procedure just described we can construct an estimator $\hat{F}_{XLR}$ for the unknown function $F_{XLR}$ allowing us to estimate how the censoring mechanism operates.
\section{Kernel density estimation in 1 and 2 dimensions}

Consider a univariate random sample $Z_1,\,Z_2, \ldots, Z_N$ from some unknown pdf $f_Z$ and suppose that the corresponding observations are all exact. The kernel density estimate $\hat{f}_Z$ of $f_Z$ is defined as
\[ \hat{f}_Z(z) = \frac{1}{Nh}\,\sum_n K\left( \frac{z-z_n}{h} \right)\]
where $h$ is an appropriately chosen parameter. The rationale for such a construction is to place a ''bump'' of size $1/N$ centered over each one of the sample values $z_n$. The general shape of each bump is determined by the choice of the kernel function $K$ while its spread is controlled by the parameter $h$. All the bumps are set to be of equal size $1/N$ due to the i.i.d. nature of the observations. The size of the bump over $z_n$ can be also interpreted as the amount of probability assigned over the interval $(z_{n-1},\,z_n]$ by the empirical cdf $\hat{F}_z$ and is thus equal to $\hat{F}_Z(z_n)-\hat{F}_Z(z_{n-1})$.

We apply the reasoning from above to the case of a univariate random sample $X_1,\,X_2, \ldots, X_N$
from some unknown density function $f_X$ such that each $X_n$ is censored into some interval $(L_n,\,R_n]$. The set $G$ is reduced to the set $\left\{\tilde{X}_1,\,\tilde{X}_2,\ldots,\tilde{X}_I \right\}$ of unique element values of the set
\[ \left\{ L_1,\, L_2, \dots, L_N \right\} \cup \left\{ R_1,\, R_2, \dots, R_N \right\} \]
listed in increasing order. We proceed to define the function $w: G\rightarrow [0,1]$ by $w(\tilde{X}_1)= \hat{F}(\tilde{X}_1)$ and $w(\tilde{X}_i)=\hat{F}(\tilde{X}_i) - \hat{F}(\tilde{X}_{i-1})$ if $2\leq i \leq I$. Now we define the smoothed density estimator $\hat{f}_X$ as
\[ \hat{f}_X(x) = \frac{1}{h}\,\sum_i w(\tilde{X}_i)\, K\left( \frac{x-\tilde{X}_i}{h} \right)\]

Next we generalize the latter construction to the case of a random sample $\{(X_n, Y_n)\}$ of censored $2$-dimensional random vectors with unknown p.d.f $f_{XY}$. The set $G$ is given by $G = G^{(x)} \times G^{(y)}$ where
\begin{eqnarray*}
G^{(x)} &=& \{\, \tilde{X}_1, \tilde{X}_2, \ldots, \tilde{X}_I \,\}\\
G^{(y)} &=& \{\, \tilde{Y}_1, \tilde{Y}_2, \ldots, \tilde{Y}_J \,\}
\end{eqnarray*}
The definition of the function $w: G\rightarrow [0,1]$ is extended as follows: $w(\tilde{X}_i,\tilde{Y}_j)=0$ if $i=1$ or $j=1$. In all other cases $w(\tilde{X}_i,\tilde{Y}_j)$ equals the cummulative probability assigned by $\hat{F}$ over the interior of the rectangle in $R^2$ defined by the points $(\tilde{X}_{i-1},\tilde{Y}_{j-1})$, $(\tilde{X}_i,\tilde{Y}_{j-1})$, $(\tilde{X}_i,\tilde{Y}_j)$ and $(\tilde{X}_{i-1},\tilde{Y}_j)$ along with the line segments connecting $(\tilde{X}_i,\tilde{Y}_{j-1})$ with $(\tilde{X}_i,\tilde{Y}_j)$ and $(\tilde{X}_{i-1},\tilde{Y}_j)$ with $(\tilde{X}_i,\tilde{Y}_j)$. Consequently, the function value $w(\tilde{X}_i,\tilde{Y}_j)$, $2\leq i \leq I$, $2\leq j \leq J$,   is given by
\[w(\tilde{X}_i,\tilde{Y}_j) = \hat{F}(\tilde{X}_i,\tilde{Y}_j) - \hat{F}(\tilde{X}_i,\tilde{Y}_{j-1}) - \hat{F}(\tilde{X}_{i-1},\tilde{Y}_j) + \hat{F}(\tilde{X}_{i-1},\tilde{Y}_{j-1})  \]
Pseudocode employing the recursive relationship from above to compute the weights $w(\tilde{X}_i,\tilde{Y}_j)$ is provided next:\\
\\
$\mbox{FOR } j=1:J$\\
$w(\tilde{X}_1,\tilde{Y}_j) = 0$\\
$\mbox{NEXT } j$\\
\\
$\mbox{FOR } i=1:I$\\
$w(\tilde{X}_i,\tilde{Y}_1) = 0$\\
$\mbox{NEXT } i$\\
\\
$\mbox{FOR } j=2:J$\\
$\mbox{FOR } i=2:I$\\
$w(\tilde{X}_i,\tilde{Y}_j) = \hat{F}(\tilde{X}_i,\tilde{Y}_j) - \hat{F}(\tilde{X}_i,\tilde{Y}_{j-1}) - \hat{F}(\tilde{X}_{i-1},\tilde{Y}_j) + \hat{F}(\tilde{X}_{i-1},\tilde{Y}_{j-1})$\\
$\mbox{NEXT } i$\\
$\mbox{NEXT } j$\\

Having developed a method for computing the weights $w(\tilde{x}_i,\tilde{y}_j)$ we are ready to present the expression for the smoothed density estimator $\hat{f}_{XY}(x, y)$:
\[ \hat{f}_{XY}(x, y) = \left(\frac{1}{h_x}\right) \left(\frac{1}{h_y}\right) \sum_{(i,j)} w(\tilde{X}_i,\,\tilde{Y}_j)\, K\left( \frac{x-\tilde{X}_i}{h_x} \right) K\left( \frac{y-\tilde{Y}_j}{h_y} \right)\]
\section{Kernel method in $M$ dimensions}
We will use $\tilde{\bold{X}}$ to designate an arbitrary element $\left(\tilde{X}_{i_1},\, \tilde{X}_{i_2}, \ldots, \, \tilde{X}_{i_M} \right)$ of the grid $G$. Furthermore, given any vector $\tilde{\bold{X}} \in G$ such that $i_m\geq 2$ for $\forall m$ we will define the vector
\[ \tilde{\bold{X}}' = \left(\tilde{X}_{i_1-1},\, \tilde{X}_{i_2-1}, \ldots, \, \tilde{X}_{i_M-1} \right)\in G \]

Let $\Omega(\bold{x})$ be the set of all hyperplanes passing through $\bold{x}$ and parallel to the coordinate planes. For example, $\Omega(\tilde{\bold{X}})$ is the set of all hyperplanes passing through $\tilde{\bold{X}}$ and parallel to the coordinate planes. Define the function $w: G\rightarrow [0,1]$ as follows: $w( \tilde{\bold{X}} ) = \hat{F}( \tilde{\bold{X}} ) = 0$ if there exists a component $\tilde{X}_{i_m}$ of $\tilde{\bold{X}}$ such that $i_m=1$. In the case when $i_m\geq 2$ for $\forall m$ the value of $w( \tilde{\bold{X}} )$ is given by the cummulative probability assigned by $\hat{F}$ over the hyperrectangle in $R^M$ bounded by the hyperplanes in $\Omega(\tilde{\bold{X}}')$ and $\Omega(\tilde{\bold{X}})$ but excluding the points lying on the hyperplanes in $\Omega(\tilde{\bold{X}}')$.

Let $\bold{h} = (h_1,\,h_2,\ldots,h_M)$. The smoothed function estimator is given by
\[ \hat{f}(\bold{x}) = \left(\prod_m \frac{1}{h_m}\right) \sum_{\tilde{\bold{x}}} w(\tilde{\bold{x}}) K\left(\bold{x};\tilde{\bold{x}},\bold{h}\right) \]
where
\[ K\left(\bold{x};\tilde{\bold{x}},\bold{h}\right) = \prod_m K\left(\frac{x_m - \tilde{x}_{i_m}}{h_m} \right) \]
\section{A loss function for computing the optimal bandwidth}

Consider once again a univariate random sample $Z_1,\,Z_2, \ldots, Z_N$ from some unknown pdf $f_Z$ and the kernel density estimator
\[ \hat{f}_Z(z) = \hat{f}_Z(z; h) = \frac{1}{N h}\,\sum_n K\left( \frac{z-Z_n}{h} \right)\]
The bandwidth $h$ will be treated as a variable for the remainder of the section. Also, to simplify notation whenever no ambiguity arises we will distinguish density functions by their argument only and drop the subscripting random variable. For example, $f(z)$ will represent $f_Z(z)$. In addition, we will use a subscript ``$-n$'' to indicate that a quantity has been derived based on the subset of the original random sample obtained after removing the $n$-th observation. For example, $\hat{f}_{-n}(z)$ is the kernel density estimator for $f(z)$ calculated after removing $Z_n$ from the original sample.

The integrated square error is defined as
\[ \int \left[\hat{f}(z)-f(z) \right]^2dz = \int \hat{f}(z)^2dz - 2\int\hat{f}(z)f(z)dz + \int f(z)^2dz \]
and the value of $h$ minimizing the risk function $R(h)$ given by
\begin{eqnarray*}
R(h) &=& E\left\{ \int \left[\hat{f}(z)-f(z) \right]^2dz \right\} \\
R(h) &=& E\left\{\int \hat{f}(z)^2dz \right\} - 2\,E\left\{\int\hat{f}(z)f(z)dz \right\} + \int f(z)^2dz
\end{eqnarray*}
is generally viewed as the optimal choice for the value of $h$ in $\hat{f}_Z(z; h)$. The term $\int f(z)^2dz$ is independent of $h$ and as a result we need to minimize 
\[ E\left\{  \int \hat{f}(z)^2dz  \right\} - 2\,E\left\{ \int\hat{f}(z)f(z)dz \right\} \]
The latter goal, however, is unachievable since the density $f_Z$ is unknown.

In reality we seek to minimize the score function
\[ M_0(h) = \int \hat{f}(z)^2dz - \frac{2}{N}\sum_n f_{-n}(Z_n) \]
for two reasons. It is straightforward to demonstrate that
\[ E\left\{ \frac{1}{n}\sum_n f_{-n}(Z_n) \right\} = E\left\{ \int\hat{f}(z)f(z)dz \right\} \]
which immediately implies that
\[E\left\{ M_0(h) \right\} = E\left\{  \int \hat{f}(z)^2dz  \right\} - 2\,E\left\{ \int\hat{f}(z)f(z)dz \right\} \]
In addition, as stated by Silverman \cite{r2} ``Assuming that the minimizer of $M_0(h)$ is close to the minimizer of $E\{M_0(h)\}$ indicates why we might hope that minimizing $M_0$ gives a good choice of smoothing parameter.''

Now we move on to motivate and introduce a score function $\tilde{M}_0(h)$ that mimics the form of $M_0(h)$ and can be used in the presence of censoring. We begin by defining the random variables
\begin{eqnarray*}
\tilde{L}_n &=& \mbox{max}\{-\infty,\,L_n\}\\
\tilde{R}_n &=& \mbox{min}\{R_n,\,+\infty\}\\
V_n &=& \frac{1}{2}(\tilde{L}_n+\tilde{R}_n)
\end{eqnarray*}
If we make the assumption that the probability distribution functions $g$ of $V_n$ and $f$ of $X_n$ are approximately equal, i.e $g(v) \approx f(v)$, then we have that
\begin{eqnarray*}
E\left\{\frac{1}{N} \sum \hat{f}_{-n}(V_n) \right\} &=& \frac{1}{N}(N) E\left\{ \hat{f}_{-1}(V_1) \right\}\\
\\
&=& E\left\{ \hat{f}_{-1}(V_1) \right\}\\
\\
&=& E\left\{ \int\hat{f}_{-1}(v)g(v)dv \right\}\\
\\
&\approx& E\left\{ \int\hat{f}_{-1}(v)f(v)dv \right\}\\
\\
&=& E\left\{ \int\hat{f}_{-1}(x)f(x)dx \right\}\\
\end{eqnarray*}
Since the expected values $E\left\{ \int\hat{f}_{-1}(x)f_X(x)dx \right\}$ and $E\left\{ \int\hat{f}(x)f_X(x)dx \right\}$ converge asymptotically we can conclude that for large samples
\[ E\left\{\frac{1}{N} \sum \hat{f}_{-n}(V_n) \right\} \approx E\left\{ \int\hat{f}_{-1}(x)f(x)dx \right\} \approx E\left\{ \int\hat{f}(x)f(x)dx \right\} \]
Consequently, we define $\tilde{M}_0(h)$ as
\[ \tilde{M}_0(h) = \int \hat{f}(x)^2dx - \frac{2}{N} \sum \hat{f}_{-n}(V_n)\]

In $M\geq 2$ dimensions we define the random variables
\begin{eqnarray*}
\tilde{L}_{nm} &=& \mbox{max}\{-\infty,\,L_{nm}\}\\
\tilde{R}_{nm} &=& \mbox{min}\{R_{nm},\,+\infty\}\\
V_{nm} &=& \frac{1}{2}(\tilde{L}_{nm}+\tilde{R}_{nm})
\end{eqnarray*}
and the random vector $\bold{V}_n = (V_{n1},\, V_{n2}, \ldots, V_{nM})$. Under the assumption that the probability distribution functions $g$ of $\bold{V}_n$ and $f$ of $\bold{X}_n$ are approximately equal, i.e $g(\bold{v}) \approx f(\bold{v})$, and based on identical reasoning we generalize the definition of $\tilde{M}_0(\bold{h})$ as follows:
\[ \tilde{M}_0(\bold{h}) = \int \hat{f}(\bold{x})^2d\bold{x} - \frac{2}{N} \sum \hat{f}_{-n}(\bold{V}_n)\]

\section{Nadaraya Watson regression with censored data}

In regression analysis the goal is to estimate the expected value $E\left\{Y|\bold{X}=\bold{x}\right\}$ based on a random sample $\left\{(\bold{X}_n,\,Y_n) \right\}$ from some unknown p.d.f. $f$ where $\bold{X}_n$ is an $M$-dimensional vector of explanatory variables. Nadaraya and Watson \cite{r3, r4} have proposed a non-parametric estimator for $E\left\{Y|\bold{X}=\bold{x}\right\}$ derived from the kernel density estimator for $f$ in the case when all sample observations are exact. We employ the newly developed censoring kernel density estimator
\[ \hat{f}(\bold{x}) = \left(\prod_m \frac{1}{h_m}\right) \sum_{\tilde{\bold{x}}} w(\tilde{\bold{x}}) K\left(\bold{x};\tilde{\bold{x}},\bold{h}\right) \]
and an identical pattern of reasoning to adapt the Nadaraya-Watson estimator for use with censored data.

In $1+1$ dimensions the censoring kernel density estimator can be written as
\[ \hat{f}(x,y) = \sum_{(i,j)} w(\tilde{x}_i, \tilde{y}_j)\frac{1}{h_x\, h_y} K\left(\frac{x-\tilde{x}_i}{h_x}\right)K\left(\frac{y-\tilde{y}_j}{h_y}\right) \]
Consequently
\begin{eqnarray*}
\hat{f}(x) &=& \int \hat{f}(x,y) dy\\
\\
&=& \sum_{(i,j)} w(\tilde{x}_i, \tilde{y}_j)\frac{1}{h_x\, h_y} K\left(\frac{x-\tilde{x}_i}{h_x}\right) \int K\left(\frac{y-\tilde{y}_j}{h_y}\right)dy\\
\\
&=& \sum_{(i,j)} w(\tilde{x}_i, \tilde{y}_j)\frac{1}{h_x\, h_y} K\left(\frac{x-\tilde{x}_i}{h_x}\right)h_y\\
\\
&=& \frac{1}{h_x}\, \sum_{(i,j)} w(\tilde{x}_i, \tilde{y}_j) K\left(\frac{x-\tilde{x}_i}{h_x}\right)
\end{eqnarray*}
and
\begin{eqnarray*}
\int y\hat{f}(x,y) dy &=& \sum_{(i,j)} w(\tilde{x}_i, \tilde{y}_j)\frac{1}{h_x} K\left(\frac{x-\tilde{x}_i}{h_x}\right) \int y\, \frac{1}{h_y} K\left(\frac{y-\tilde{y}_j}{h_y}\right)dy\\
\\
&=& \frac{1}{h_x}\,\sum_{(i,j)} w(\tilde{x}_i, \tilde{y}_j) K\left(\frac{x-\tilde{x}_i}{h_x}\right)\tilde{y}_j
\end{eqnarray*}
Now we define the estimator $E\left\{Y|X=x\right\}$ as follows:
\[ E\left\{Y|X=x\right\} = \frac{\int y\hat{f}(x,y) dy}{\hat{f}(x)} = \frac{\sum_{(i,j)} w(\tilde{x}_i, \tilde{y}_j) K\left(\frac{x-\tilde{x}_i}{h_x}\right)\tilde{y}_j}{\sum_{(i,j)} w(\tilde{x}_i, \tilde{y}_j) K\left(\frac{x-\tilde{x}_i}{h_x}\right)} \]

In $(M+1)$ dimensions the same reasoning leads us to define the estimator $E\left\{Y|\bold{X}=\bold{x}\right\}$ as
\[ E\left\{Y|\bold{X}=\bold{x}\right\} = \frac{\int y\hat{f}(\bold{x},y) dy}{\hat{f}(\bold{x})} = \frac{\sum_{(\tilde{\bold{x}},y_j)} w(\tilde{\bold{x}}, \tilde{y}_j) K\left(\bold{x};\tilde{\bold{x}},\bold{h}\right)\tilde{y}_j}  {\sum_{(\tilde{\bold{x}},y_j)} w(\tilde{\bold{x}}, \tilde{y}_j) K\left(\bold{x};\tilde{\bold{x}},\bold{h}\right)} \]

\section{Parameter estimation for a multinomial distribution in the presence of censoring}

Let $c_1$ and $c_2$ be the respective observed numbers of outcomes of type 1 and type 2 in a binomial experiment with $N$ trials, $u = N-c_1-c_2\geq 1$ number of trials with unknown outcomes and probability $\pi$ of a single trial being of type 1. Let $N_1$ and $N_2$ designate the actual counts of type 1 and type 2. Consequently $N_1$ and $N_2$ are censored such that $(N_1,\,N_2)\in S_2$ where the set $S_2$ is defined by
\[ S_2 = \{(l_1, l_2) \,|\, l_1,l_2 \mbox{ are non-negative integers},\,  l_1\geq c_1,\, l_2\geq c_2,\, l_1+l_2=N\} \]
If $E$ designates the event $(N_1,\,N_2)\in S$ then the likelihood of observing $E$ is given by
\[L(\pi; E) = \sum_{(n_1,n_2)\in S_2}  \frac{N!}{n_1!\,n_2!}\,(\pi)^{n_1}\,(1-\pi)^{n_2} = \sum^{c_1+u}_{n_1=c_1}  {N \choose n_1}(\pi)^{n_1} (1-\pi)^{N-n_1}\]
As already derived, the value $\hat{\pi}$ of $p$ that maximizes the function
\[L(p; E) = \sum^{c_1+u}_{n_1=c_1}  {N \choose n_1}(p)^{n_1} (1-p)^{N-n_1}\]
subject to the constraint $0 \leq p \leq 1$ is given by
\[ \hat{\pi} = \left\{
	\begin{array}{lll}
		0 & \mbox{if} & c_1=0 \mbox{ and } c_2\geq 1 \\ 
		\\
		1 & \mbox{if} & c_2=0 \mbox{ and } c_1\geq 1 \\
		\\
		1/2 & \mbox{if} & u=N \\
		\\		\left(1+\sqrt[u+1]{\frac{c_2 (c_2+1)\ldots(c_2+u)}{c_1 (c_1+1)\ldots(c_1+u)}}\,\right)^{-1}&\mbox{o.w.}
	\end{array}
\right. \]

The treatment of an multinomial experiment with $N$ trials, $M$ possible outcome types and probabilities $\pi_1, \pi_2, \ldots, \pi_M$ of each outcome type is based on the same reasoning. We use $c_1, c_2, \ldots, c_M$ to designate the observed counts of each type and $N_1, N_2, \ldots, N_M$ to designate the actual and possibly censored outcome counts. Suppose $u = N - \sum_m c_m \geq 1$ and define the vectors
\begin{eqnarray*}
\bold{p} &=& (p_1,\,p_2,\ldots,p_M)\\
\bold{c} &=& (c_1,\,c_2,\ldots,c_M)\\
\bold{n} &=& (n_1,\,n_2,\ldots,n_M)\\
\bold{N} &=& (N_1,\,N_2,\ldots,N_M)
\end{eqnarray*}
The definition of the set $S_2$ generalizes to
\[ S_M = \{(l_1, l_2, \ldots, l_M) \,|\, l_{m}\mbox{ is a non-negative integer},\,  l_{m}\geq c_{m},\, \sum_m l_m=N \} \]
and accordingly $E$ is redefined to be the event $(N_1,\,N_2, \ldots, N_M)\in S_M$. The likelihood of $E$ as a function of $\bold{p}$ is given by
\[L(\bold{p}; E) = \sum_{\bold{n}\in S_M}  \frac{N!}{n_1!\,n_2!\ldots n_M!}\,(p_1)^{n_1}\,(p_2)^{n_2}\ldots \,(p_M)^{n_M}\]
An approximate solution to the resulting estimation problem can be constructed as follows. If $\hat{p}_m$ is the value of the variable $p_m$ that maximizes the function
\[ \sum^{c_m+u}_{n_m=c_m}  {N \choose n_m}(p_m)^{n_m} (1-p_m)^{N-n_m} \]
then we could employ
\[\hat{\pi}^*_m = \frac{\hat{p}_m}{\hat{p}_1+\hat{p}_2+\ldots+\hat{p}_M}\]
as an estimator for the unknown probability $\pi_m$.

Next we consider a trinomial $(M=3)$ experiment such that $u_{12}$ trials are of type $1$ or type $2$ and $u_{23}$ are of type $2$ or type $3$ and define
\begin{eqnarray*}
u_1 &=& u_{12}\\
u_2 &=& \mbox{min}\{ N - c_1 -c_2 - c_3,\, u_{12}+u_{23} \}\\
u_3 &=& u_{23}
\end{eqnarray*}
Let $\hat{p}_m$ be the value of the variable $p_m$ that maximizes the function
\[ \sum^{c_m+u_m}_{n_m=c_m}  {N \choose n_m}(p_m)^{n_m} (1-p_m)^{N-n_m} \]
and
\[\hat{\pi}^*_m = \frac{\hat{p}_m}{\hat{p}_1+\hat{p}_2+\ldots+\hat{p}_M}\]
The quantities $\hat{\pi}^*_1$, $\hat{\pi}^*_2$ and $\hat{\pi}^*_3$ can be used to estimate the unknown probabilities $\pi_1$, $\pi_2$ and $\pi_3$. Generalizing to the case of $M$ possible outcomes in the presence of partial censoring is straightforward. Let $u_m$ be the maximum possible number of censored outcomes of type $m$ and assume that $1, 2, \ldots, u_m$ are all possible counts for the number of unobserved outcomes of type $m$. Consequently $\hat{\pi}^*_m$ is a potential estimator for $\pi_m$.

So far we have been constructing likelihood functions without making assumptions or having the benefit of prior knowledge about the nature of the censoring mechanism. Let $q_m$ be the conditional probability of observing an outcome ot type $m$ and $\bold{q}=(q_1,\,q_2,\ldots,q_M)$. For example, let us consider a binomial ($M=2$) experiment with known parameters $q_1$ and $q_2$. The probability of not being able to observe the outcome of a single trial $X_n$ is given by $(1-q_1)p_1+(1-q_2)p_2 = (p_1+p_2)-p_1q_1-p_2q_2$. Consequently the likelihood of observing $c_1$ outcomes of type $1$, $c_2$ outcomes of type $2$ and $u = N-c_1-c_2$ outcomes of unknown type is
\[L(\bold{p}, \bold{q};\bold{c},\,u) = \frac{N!}{c_1!\,c_2!\,u!}\,(p_1\,q_1)^{c_1}\,(p_2q_2)^{c_2}\left[(1-q_1)\,p_1+(1-q_2)\,p_2\right]^u\]
Generalizing is trivial:
\[L(\bold{p}, \bold{q};\bold{c},\,u) = \frac{N!}{c_1!\,c_2! \ldots c_M!\,u!}\,\prod_m{(p_m\,q_m)^{c_1}}\,\left[\sum_m{(1-q_m)\,p_m}\right]^u\]
where $u = N - \sum_m c_m$.

Finally we turn our attention to a binomial experiment such that $q_1$ remains unknown but $q_2$ is known. The outcome $x_n$ of a single trial $X_n$ can be classified in exactly one of the following four categories: observed of type 1, observed of type 2, unobserved of type 1 and unobserved of type 2. Let $\tilde{N}_1$ designate the number of censored outcomes of type 1, $\tilde{N}_2$ designate the number of censored outcomes of type 2 and the set $\tilde{S}_2$ be defined as
\[ \tilde{S}_2 = \{(l_1, l_2) \,|\, l_1,\,l_2\mbox{ are non-negative integers},\,  l_1+l_2=N-c_1-c_2 \} \]
The likelihood of the event $\tilde{E}=\,$``$(\tilde{N}_1,\,N_2)\in \tilde{S}_2$'' is given by
\[ L(\bold{\pi}, \bold{q}; \tilde{E}) = \sum_{(\tilde{n}_1,\tilde{n}_2)}  \frac{N!}{c_1!\,c_2!\,\tilde{n}_1!\,\tilde{n}_2!}\,(\pi_1q_1)^{c_1}\,(\pi_2q_2)^{c_2}\left[(1-q_1)\pi_1\right]^{\tilde{n}_1}\left[(1-q_2)\pi_2\right]^{\tilde{n}_2} \]
where the summation index $(\tilde{n}_1,\tilde{n}_2)$ spans the set $\tilde{S}_2$. Consequently we seek to maximize the function
\[ L(\bold{p},\,q'_2; E) = \sum_{(\tilde{n}_1,\tilde{n}_2)}  \frac{N!}{c_1!\,c_2!\,\tilde{n}_1!\,\tilde{n}_2!}\,(p_1q_1)^{c_1}\,(p_2q_2')^{c_2}\left[(1-q_1)p_1\right]^{\tilde{n}_1}\left[(1-q_2')p_2\right]^{\tilde{n}_2} \]
subject to the constraints $p_1+p_2=1$ and $0\leq q_2' \leq 1$.

\section{Analysis of contingency tables with incomplete counts}

Since each cell in an $I\times J$ contingency table can be uniquely associated with an ordered pair $(i,j)$ the set of ordered pairs $\{(i,j)\}$ constitutes the space of possible outcomes for a sample random variable $X_n$. Define the probabilities $\pi_{ij}$, $q_{ij}$ and $\alpha_{ij}$ as 
\begin{eqnarray*}
\pi_{ij} &=& \mbox{Prob}\left\{ X_n = (i,j) \right\}\\
q_{ij} &=&  \mbox{Prob}\left\{ X_n \mbox{ is observed}\,|\,X_n=(i,j) \right\}\\
\alpha_{ij} &=& \mbox{Prob}\left\{ X_n=(i,j) \mbox{ and }X_n \mbox{ is observed} \right\} = \pi_{ij} q_{ij}
\end{eqnarray*}
Furthermore, let $c_{ij}$ and $N_{ij}$ be the respective observed and actual counts in cell $(i,j)$. We can quantify the effect of the censoring mechanism by observing that the ratio $\hat{\alpha}_{ij} = \frac{c_{ij}}{N}$ constitutes an MLE for the joint probability $\alpha_{ij}$ and using the plug-in principle within the equation $\alpha_{ij} = \pi_{ij} q_{ij}$ to obtain the estimator $\hat{q}_{ij} = \frac{c_{ij}}{\hat{\pi}_{ij}N}$ for the unknown probability $q_{ij}$.

The actual count $N_{ij}$ may be unknown due to the censoring mechanism. From the definitions follows that $c_{ij} = N_{ij}$ if outcomes of type $(i,j)$ are not subject to censoring and $c_{ij} \leq N_{ij}$ otherwise. Finally, let us use $N_j = \sum_i N_{ij}$ to designate the $j$-th column total and in the case when $N_j$ is known let $u_j = N_j - \sum_i c_{ij}$ designate the number of sample outcomes censored into the $j$-th column. 

Consider the special case of a $2\times 2$ ($I=2, J=2$) contingency table and the null hypothesis
\[ H_0:\,\mbox{Prob}\left\{ X_n = (1,1)\,|\, X_n \in \{(1,1),(2,1)\} \right\} = \mbox{Prob}\left\{ X_n = (1,2)\,|\, X_n \in \{(1,2),(2,2)\} \right\} \]
which can be rewritten as
\[ H_0:\, \frac{\pi_{11}}{\pi_{11}+\pi_{21}} = \frac{\pi_{12}}{\pi_{12}+\pi_{22}} \]
Assuming $H_0$ in an estimation problem amounts to introducing the constraint
\[ \frac{p_{11}}{p_{11}+p_{21}} = \frac{p_{12}}{p_{12}+p_{22}} \]
where $p_{ij}$ is the variable associated with the unknown cell probability $\pi_{ij}$. In the special case of predetermined column totals $N_1$ and $N_2$ we have that $\pi_{11}+\pi_{12}=1$ as well as $\pi_{12}+\pi_{22}=1$. Consequently, the null hypothesis is reduced to $H_0:\, \pi_{11}=\pi_{12}=\pi$ and accordingly the null constraint becomes $p_{11}=p_{22}=p$.

Before turning our attention to three examples of censored $2\times 2$ contingency tables we introduce some additional notation:
\begin{eqnarray*}
&S& = \{(l_{11},\,l_{21},\,l_{12},\,l_{22}) \,|\, l_{ij}\mbox{ is a non-negative integer},\,  l_{ij}\geq c_{ij},\, \sum_{(i,j)} l_{ij}=N \}\\
&\bold{N}& = (N_{11},\,N_{21},\,N_{12},\,N_{22})\\
&\bold{n}& = (n_{11},\,n_{21},\,n_{12},\,n_{22})\\
&\bold{c}& = (c_{11},\,c_{21},\,c_{12},\,c_{22})\\
&\bold{p}& = (p_{11},\,p_{21},\,p_{12},\,p_{22})\\
&\bar{\bold{c}}& = (\bar{c}_{11},\,\bar{c}_{21},\,\bar{c}_{12},\,\bar{c}_{22})
\end{eqnarray*}
In each example we construct the likelihood necessary to derive a set of estimators $\{\hat{\pi}_{ij}\}$ for the elements of $\{\pi_{ij}\}$. A superscript ``$(0)$'' will be used to label quantities derived under $H_0$. For example, $\hat{\pi}^{(0)}_{ij}$ is the null esimator for $\pi_{ij}$.

\subsection{Example 1}

Suppose that $N_1$ and $N_2$ are predetermined by the experimenter, the counts $N_{11}$ and $N_{21}$ are exact implying $u_1=0$ while the counts $N_{12}$ are $N_{22}$ are censored implying $u_2\geq 1$. Let $E_1$ designate the event ``$\bold{N}\in S$ and $N_{11}+N_{21}=N_1$ and $N_{12}+N_{22}=N_2$''. The likelihood of observing $E_1$ is given by
\[ L(\bold{p}; E_1) = L_1(p_{11})\,L_2(p_{12})\]
where
\begin{eqnarray*}
L_1(p_{11}) &=& {N_1 \choose N_{11}}(p_{11})^{N_{11}} (1-p_{11})^{N_1-N_{11}} \\
\\
L_1(p_{12}) &=& \sum^{c_{12}+u_2}_{n_{12}=c_{12}}  {N_2 \choose n_{12}}(p_{12})^{n_{12}} (1-p_{12})^{N_2-n_{12}}
\end{eqnarray*}
Since the column totals $N_1$ and $N_2$ are fixed and known in advance, under $H_0$ the likelihood function needs to be modified by setting $p = p_{11} = p_{12}$:
\begin{eqnarray*}
L(\bold{p}; E, H_0) &=& \sum^{c_{12}+u_2}_{n_{12}=c_{12}}  {N_1 \choose N_{11}}(p)^{N_{11}} (1-p)^{N_1-N_{11}} \, {N_2 \choose n_{12}}(p)^{n_{12}} (1-p)^{N_2-n_{12}}  \\
\\
L(\bold{p}; E, H_0) &=& \sum^{c_{12}+u_2}_{n_{12}=c_{12}} {N_1 \choose N_{11}}{N_2 \choose n_{12}}(p)^{N_{11}+n_{12}} (1-p)^{N-N_{11}-n_{12}}
\end{eqnarray*}

Let $\bold{t} = (\bar{N}_{11},\,\bar{c}_{12},\,\bar{c}_{22},\,\bar{u}_2)$ be a particular vector of counts for the contingency table. Then the probability of observing $\bold{t}$ is given by
\[ \mbox{Prob}\{\bold{t}\} = P_1(\bar{N}_{11})\,P_2(\bar{c}_{12},\,\bar{c}_{22},\,\bar{u}_2) \]
where
\begin{eqnarray*}
P_1(\bar{N}_{11}) &=& \frac{N_1}{\bar{N}_{11}!\,\bar{N}_{21}!}\,(\pi_{11})^{\bar{N}_{11}}\,(\pi_{21})^{\bar{N}_{21}}\\
\\
P_2(\bar{c}_{12},\,\bar{c}_{22},\,\bar{u}_2) &=& \frac{N_2}{\bar{c}_{12}!\,\bar{c}_{22}!\,\bar{u}_2!}\,(\alpha_{12})^{\bar{c}_{12}}\,(\alpha_{22})^{\bar{c}_{22}}\,(1-\alpha_{12}-\alpha_{22})^{\bar{u}_2}
\end{eqnarray*}
We can estimate $\mbox{Prob}\{\bold{t}\}$ by using $\hat{\pi}_{11},\, \hat{\pi}_{21},\, \hat{\alpha}_{12}$ and $\hat{\alpha}_{22}$ for the unknown probabilities $\pi_{11},\, \pi_{21},\, \alpha_{12}$ and $\alpha_{22}$. Under $H_0$ we estimate $\mbox{Prob}\{\bold{t}\}$ by employing the appropriate null estimators $\hat{\pi}^{(0)}_{11}$ and $\hat{\pi}^{(0)}_{21}$ as opposed to $\hat{\pi}_{11}$ and $\hat{\pi}_{21}$.

\subsection{Example 2}
Suppose $N_1$ and $N_2$ are predetermined by the experimenter and the counts $N_{11}, N_{21}, N_{12}, N_{22}$ are all unobserved. We use $E_2$ designate the event ``$\bold{N}\in S$ and $N_{11}+N_{21}=N_1$ and $N_{12}+N_{22}=N_2$''. The likelihood of $E_2$ is given by
\[ L(\bold{p}; E_2) = L_1(p_{11})\,L_2(p_{12})\]
where
\begin{eqnarray*}
L_1(p_{11}) &=& \sum^{c_{11}+u_1}_{n_{11}=c_{11}}  {N_1 \choose n_{11}}(p_{11})^{n_{11}} (1-p_{11})^{N_1-n_{11}} \\
\\
L_1(p_{12}) &=& \sum^{c_{12}+u_2}_{n_{12}=c_{12}}  {N_2 \choose n_{12}}(p_{12})^{n_{12}} (1-p_{12})^{N_2-n_{12}}
\end{eqnarray*}
The two factors $L_1(p_{11}; E)$ and $L_2(p_{12}; E)$ can be maximized independently if no further assumptions are made. Under the null constraint $p = p_{11} = p_{12}$ the likelihood $L(\bold{p}; E_2)$ is modified as follows:
\begin{eqnarray*}
L(\bold{p}; E_2, H_0) &=& \sum_{(n_{11},\,n_{12})} {N_1 \choose n_{11}}(p)^{n_{11}} (1-p)^{N_1-n_{11}} \, {N_2 \choose n_{12}}(p)^{n_{12}} (1-p)^{N_2-n_{12}} \\
\\
L(\bold{p}; E_2, H_0) &=& \sum_{(n_{11},\,n_{12})} {N_1 \choose n_{11}}{N_2 \choose n_{12}}(p)^{n_{11}+n_{12}}\, (1-p)^{N-n_{11}-n_{12}}
\end{eqnarray*}
where $(n_{11},\,n_{12}) \in S$ and $n_{11}+n_{21}=N_1$ and $n_{12}+n_{22}=N_2$.

The probability of a particular contingency table configuration is given by
\[ \mbox{Prob}\{\bar{\bold{c}},\,\bar{u}_1,\,\bar{u}_2 \} = P_1(\bar{c}_{11},\,\bar{c}_{21},\,\bar{u}_1)\,P_2(\bar{c}_{12},\,\bar{c}_{22},\,\bar{u}_2) \]
with
\begin{eqnarray*}
P_1(\bar{c}_{11},\,\bar{c}_{21},\,\bar{u}_1) &=& \frac{N_1}{\bar{c}_{11}!\,\bar{c}_{21}!\,\bar{u}_1!}\,(\alpha_{11})^{\bar{c}_{11}}\,(\alpha_{21})^{\bar{c}_{21}}\,(1-\alpha_{11}-\alpha_{21})^{\bar{u}_1}\\
\\
P_2(\bar{c}_{12},\,\bar{c}_{22},\,\bar{u}_2) &=& \frac{N_2}{\bar{c}_{12}!\,\bar{c}_{22}!\,\bar{u}_2!}\,(\alpha_{12})^{\bar{c}_{12}}\,(\alpha_{22})^{\bar{c}_{22}}\,(1-\alpha_{12}-\alpha_{22})^{\bar{u}_2}
\end{eqnarray*}
The estimators for $\alpha_{ij}$ remain unchanged under $H_0$ unless additional assumptions are made regarding the nature of the censoring mechanism.

\subsection{Example 3}

Suppose that $N_{11}, N_{21}, N_{12}, N_{22}$ as well as the column totals $N_1$ and $N_2$ are all unknown. Let
$u = N - (c_{11}+c_{21}+c_{12}+c_{22})$ and $E_3$ designate the event ``$\bold{N}\in S$''. The estimators $\hat{\pi}_{ij}$ maximize the likelihood
\[L(\bold{p}; E_3) = \sum_{\bold{n}\in S}  \frac{N!}{n_{11}!\,n_{21}!\,n_{12}!\,n_{22}!}\,(p_{11})^{n_{11}}\,(p_{21})^{n_{21}} \,(p_{12})^{n_{12}}\,(p_{22})^{n_{22}}\]
By enforcing the null constraint the above likelihood is reduced to
\begin{eqnarray*}
L(\bold{p}; E_3, H_0) &=& \sum_{\bold{n}\in S}  \frac{N!}{n_{11}!\,n_{21}!\,n_{12}!\,n_{22}!}\,(p)^{n_{11}}\,(1-p)^{n_{21}} \,(p)^{n_{12}}\,(1-p)^{n_{22}}
\\
L(\bold{p}; E_3, H_0) &=& \sum_{\bold{n}\in S}  \frac{N!}{n_{11}!\,n_{21}!\,n_{12}!\,n_{22}!}\,(p)^{n_{11}+n_{12}}\,(1-p)^{n_{21}+n_{22}}
\end{eqnarray*}
The probability of a particular contingency table configuration is given by
\[ \mbox{Prob}\{ \bar{\bold{c}},\,\bar{u} \} =  \frac{N}{\bar{c}_{11}!\,\bar{c}_{21}!\,\bar{c}_{12}!\,\bar{c}_{22}!\,\bar{u}!}\, P_1(\bar{\bold{c}})\,P_2(u) \]
where
\begin{eqnarray*}
P_1(\bar{\bold{c}}) &=& (\alpha_{11})^{\bar{c}_{11}}\,(\alpha_{21})^{\bar{c}_{21}}\,(\alpha_{12})^{\bar{c}_{12}}\,(\alpha_{22})^{\bar{c}_{22}}\\
\\
P_2(u) &=& (1-\alpha_{11}-\alpha_{21}-\alpha_{12}-\alpha_{22})^{\bar{u}}
\end{eqnarray*}
Assuming $H_0$ does not modify the estimate for $\mbox{Prob}\{ \bar{\bold{c}},\,\bar{u} \}$.\\

Extending the ideas presented in this section to the construction of appropriate likelihood functions for contingency tables with $I\geq 2$ rows and $J\geq 2$ columns in the presence of a censoring mechanism should be trivial in most cases. Solving the resulting optimizataion problems, however, may be far from straightforward.\\
%


\end{document}